\documentclass[sn-mathphys]{sn-jnl}

\jyear{2021}%

\usepackage{amsmath}
\usepackage{mathtools}
\usepackage{amssymb}
\usepackage{amsthm}
\usepackage{graphicx}
\usepackage{enumerate}
\usepackage[mathscr]{eucal}
\usepackage{hyperref}
\hypersetup{
    colorlinks=true,
    linkcolor=blue,
    filecolor=magenta,      
    urlcolor=cyan,
}

\numberwithin{equation}{section}
\usepackage[english]{babel}
\DeclareMathOperator{\conv}{conv}
\DeclareMathOperator{\sgn}{sgn}

\theoremstyle{thmstyleone}%
\newtheorem{theorem}{Theorem}[section]
\newtheorem{prop}[theorem]{Proposition}%
\newtheorem{cor}[theorem]{Corollary}
\newtheorem{lemma}[theorem]{Lemma}

\theoremstyle{thmstyletwo}%

\theoremstyle{thmstylethree}%
\newtheorem{definition}{Definition}%

\begin{document}

\title[Birkhoff-James Orthogonality in Sequence Spaces]{Birkhoff-James Orthogonality and Its Local Symmetry in Some Sequence Spaces}


\author*[1]{\fnm{Babhrubahan} \sur{Bose}}\email{babhrubahanb@iisc.ac.in}

\author[2]{\fnm{Saikat} \sur{Roy}}\email{saikatroy.cu@gmail.com}
\equalcont{These authors contributed equally to this work.}

\author[1]{\fnm{Debmalya} \sur{Sain}}\email{saindebmalya@gmail.com}
\equalcont{These authors contributed equally to this work.}

\affil*[1]{\orgdiv{Department of Mathematics}, \orgname{Indian Institute of Science}, \orgaddress{\street{C V Raman Road}, \city{Bangalore}, \postcode{560012}, \state{Karnataka}, \country{India}}}

\affil[2]{\orgdiv{Department of Mathematics}, \orgname{National Institute of Technology}, \orgaddress{\street{Mahatma Gandhi Road}, \city{Durgapur}, \postcode{713209}, \state{West Bengal}, \country{India}}}

\affil[3]{\orgdiv{Department of Mathematics}, \orgname{Indian Institute of Science}, \orgaddress{\street{C V Raman Road}, \city{Bangalore}, \postcode{560012}, \state{Karnataka}, \country{India}}}


\abstract{ We study Birkhoff-James orthogonality and its local symmetry in some sequence spaces namely $\ell_p,$ for $1\leq p\leq\infty$, $p\neq2$, $c$, $c_0$ and $c_{00}$. Using the characterization of the local symmetry of Birkhoff-James orthogonality, we characterize isometries of each of these spaces onto itself and obtain the Banach-Lamperti theorem for onto operators on the sequence spaces.}





\maketitle

\section{Introduction}
{In recent times, symmetry of Birkhoff-James orthogonality has been a topic of considerable interest \cite{annal},  \cite{dkp}, \cite{1}, \cite{3},  \cite{4},  \cite{5}, \cite{8}. It is now well known that the said symmetry plays an important role in the study of the geometry of Banach spaces. The present article aims to explore Birkhoff-James orthogonality and its local symmetry in some well studied sequence spaces.  As an outcome of our exploration, we acquire the Banach-Lamperti Theorem \cite{BL} for onto operators on some classical sequence spaces by characterizing the onto isometries of the same. We would like to mention that recently such a study has been carried out in the context of $\ell^n_p$ spaces for $1\leq p\leq \infty$, $p\neq2$ in \cite{CSS}.}

Let us now establish the relevant notations and terminologies to be used throughout the article. Denote the scalar field $\mathbb{R}$ or $\mathbb{C}$ by $\mathbb{K}$ and recall the sign function $\sgn:\mathbb{K}\to\mathbb{K},$ given by
\[\sgn(x)=
\begin{cases}
\frac{x}{\lvert x\rvert},\;x\neq0,\\
0,\;x=0.
\end{cases}\] 
Consider a normed linear space $\mathbb{X}$ over $\mathbb{K}$ and denote its continuous dual by $\mathbb{X}^*$. Let $J(x)$ denote the collection of all support functionals of a non-zero $x\in \mathbb{X}$, i.e.,
\begin{align}\label{support}
    J(x):=\{f\in \mathbb{X}^*:\|f\|=1,\;\lvert f(x)\rvert=\|x\|\}.
\end{align}
A non-zero {element} $x\in\mathbb{X}$ is said to be \textit{smooth} if $J(x)$ is singleton.\par
Given $x,y\in \mathbb{X}$, $x$ is said to be \textit{Birkhoff-James orthogonal} to $y$ \cite{B}, denoted by $x\perp_By$, if
\begin{align*}
    \|x+\lambda y\|\geq\|x\|,~~\textit{for~all~}\lambda\in\mathbb{K}.
\end{align*}
James proved in \cite{james} that $x\perp_By$ if and only if $x=0$ or there exists $f\in J(x)$ such that $f(y)=0$. In the same article he also proved that a non-zero $x\in \mathbb{X}$ is smooth if and only if Birkhoff-James orthogonality is right additive at $x$, i.e.,
\begin{align*}
    x\perp_By,~x\perp_Bz~~\Rightarrow~~x\perp_B(y+z),~~\textit{for every}~y,z\in\mathbb{X}.
\end{align*}
\par
Birkhoff-James orthogonality is not symmetric in general, i.e., $x\perp_By$ does not necessarily imply that $y\perp_Bx$. In fact, James proved in \cite{james2} that Birkhoff-James orthogonality is symmetric in a normed linear space of dimension higher than 2 if and only if the space is an inner product space. However, the importance of studying the local symmetry of Birkhoff-James orthogonality in describing the geometry of normed linear spaces has been illustrated in \cite[Theorem 2.11]{CSS}, \cite[Corollary 2.3.4.]{Sain}. Let us recall the following definition in this context from \cite{Sain2}, which is of paramount importance in our present study.
\begin{definition}
An element $x$ of a normed linear space $\mathbb{X}$ is said to be \textit{left-symmetric} (\textit{resp. right-symmetric}) if 
\begin{align*}
    x\perp_By\;\Rightarrow\; y\perp_Bx~~(\textit{resp.~}y\perp_Bx\;\Rightarrow\;x\perp_By),
\end{align*}
for every $y\in \mathbb{X}$.
\end{definition}
The left-symmetric and the right-symmetric points of $\ell_p^n$ spaces where $1\leq p\leq \infty$, $p\neq2,$ were characterized in \cite{CSS}. {Here we take a step forward towards generalizing these results in the following sequence spaces: $\ell_p,$ for $1\leq p\leq\infty$ and $p\neq 2$, $c$, $c_0$ and $c_{00}$.} Characterizations of the smooth points, the left-symmetric points and the right-symmetric points of a given Banach space are of paramount importance in understanding the geometry of the Banach space. We refer the readers to \cite{annal}, \cite{dkp}, \cite{1}, \cite{3}, \cite{4}, \cite{5}, \cite{8}, \cite{10}, \cite{SRBB}, \cite{12}, \cite{turnsek} for {some prominent work in this direction.} \par
The local symmetry of Birkhoff-James orthogonality in a Banach space also plays an important role in determining the isometric isomorphisms on the space. Let us observe that Corollary $2.3.4.$ of \cite{Sain} in this regard can be stated in the following generalized form:
\begin{cor}\label{isometry}
Let $\mathbb{X}$ and $\mathbb{Y}$ be two normed linear spaces and let $T:\mathbb{X}\to\mathbb{Y}$ be an onto linear isometry. Then $x\in\mathbb{X}$ is left-symmetric (resp. right-symmetric) if and only if $T(x)\in\mathbb{Y}$ is left-symmetric (resp. right-symmetric). 
\end{cor}
This result is used for proving the Banach-Lamperti {Theorem} for onto operators on the sequence spaces, i.e., for the case where the measure space is $\mathbb{N}$ equipped with the counting measure by finding the onto isometries of $\ell_p$, for $1\leq p\leq\infty$ and $p\neq 2$. We also do the same for the spaces $c$, $c_0$ and $c_{00}$ as a direct consequence of the results characterizing the local symmetry of Birkhoff-James orthogonality in these spaces. It can be noted that Lamperti's idea in \cite{BL} uses the concept of convexity, concavity and Radon-Nikodym derivatives along with the properties of the integral involved in the definition of the $L_p$ norm and therefore cannot be generalized in case of $p=\infty$. Our approach using the local symmetry of Birkhoff-James orthogonality however, has no such restrictions and hence is applied for the $p=\infty$ case as well.\par
In the first section we completely characterize Birkhoff-James orthogonality in $\ell_\infty$ over $\mathbb{K}$ and then characterize the left-symmetric and the right-symmetric points of the space. As a corollary of our results, we obtain characterizations of Birkhoff-James orthogonality and the left-symmetric and the right-symmetric points in $c$, $c_0$ and $c_{00}$. Using Corollary \ref{isometry}, we find the isometries of each of these spaces onto itself.\par
In the second and third sections we {obtain} the same characterizations for $\ell_1$ and $\ell_p$ spaces for $1<p<\infty$ and $p\neq2$ respectively. Observe that the $p=2$ case is trivial since $\ell_2$ is a Hilbert space. We also find all the isometries of these spaces onto themselves using Corollary \ref{isometry}. \par
Since we are proving the Banach-Lamperti {Theorem} for onto operators on the sequence spaces by this isometry characterization, we define signed permutation operators on $\mathbb{K}^\mathbb{N}$, the vector space of all sequences in $\mathbb{K}$.
\begin{definition}
A map $T:\mathbb{K}^\mathbb{N}\to\mathbb{K}^\mathbb{N}$ is said to be a \textit{signed permutation operator} if there exists {a bijection} $\sigma:\mathbb{N}\to\mathbb{N}$  such that
\begin{equation*}
    T(x)=\left(c_nx_{\sigma(n)}\right),~~x=\left(x_n\right)_{n\in\mathbb{N}}\in\mathbb{K}^\mathbb{N},
\end{equation*}
where $\lvert c_n\rvert =1,$ for every $n\in\mathbb{N}$.
\end{definition}

\section{Geometry of $\ell_\infty$}
In this section, we characterize Birkhoff-James orthogonality between two elements of $\ell_\infty$ and then obtain characterizations of the smooth points, the left-symmetric points and the right-symmetric points of the space. To serve our purpose, {we review some basic facts about the convergence of a $\mathbb{K}$-valued sequence under an ultrafilter on $\mathbb{N}$. A detailed treatment on ultrafilters can be found in \cite{Comfort}, \cite{Comfort N}.}

\subsection{Ultrafilters on $\mathbb{N}$ and convergence of sequences under them}\hfill
\\

We begin by recalling a few definitions.
\begin{definition}[Filters and Ultrafilters]
A non-empty subset $\mathcal{F}$ of the power set of a non-empty set $X$ is said to be a \textit{filter} on $X$ if\\
$(i)$ $\emptyset\notin\mathcal{F}.$\\
$(ii)$ $A,\,B\in\mathcal{F}~~\Rightarrow~~A\cap B\in\mathcal{F}.$\\
$(iii)$ $A\in\mathcal{F}$ \textit{and} $A\subset B~~\Rightarrow B\in\mathcal{F}$.\\
A filter $\mathcal{U}$ on $X$ is said to be an \textit{ultrafilter} on $X$ if any filter on $X$ containing $\mathcal{U}$ must be $\mathcal{U}$.
\end{definition}
Note that a filter $\mathcal{U}$ is an ultrafilter if and only if for every $A\subset X$, $A\in\mathcal{U}$ or $X\backslash A\in\mathcal{U}$ if and only if for every $n\in\mathbb{N}$ and $A_1,A_2,\dots A_n\subset X$, $X=\bigcup\limits_{i=1}^nA_i$ implies $A_i\in \mathcal{U},$ for some $1\leq i\leq n$. 
An ultrafilter $\mathcal{U}$ on $X$ is called a \textit{principal ultrafilter} if there exists $x\in X$ such that $\{x\}\in\mathcal{U}$. An ultrafilter which is not a principal ultrafilter is called a \textit{free ultrafilter}.\par 
We also recall the definition of a filter base.
\begin{definition}
A non-empty subset $\mathcal{B}$ {of the power set of a non-empty set $X$} is said to be a filter base if\\
$(i)$ $\emptyset\notin \mathcal{B}$.\\
$(ii)$ If $A,\,B\in\mathcal{B}$, there exists $C\in\mathcal{B}$ such that $C\subset A\cap B$. 
\end{definition}
Note that every filter base $\mathcal{B}$ is contained in a unique minimal filter given by $\{A\subset X:B\subset A,~\textit{for~some~}B\in\mathcal{B}\}$. Since every filter is contained in some ultrafilter by Zorn's lemma, every filter base is also contained in some ultrafilter.\par
We now focus on the case $X=\mathbb{N}$. Recall the definition of convergence along a filter:
\begin{definition}[Convergence along a filter]
Let $(x_n)_{n\in\mathbb{N}}$ be a sequence in $\mathbb{K}$ and let $\mathcal{F}$ be a filter on $\mathbb{N}$. Then we say $x_n$ converges to some point $x_0\in\mathbb{K}$, denoted by $\lim\limits_{\mathcal{F}}x_n$, under $\mathcal{F}$ if for every $\epsilon>0$, 
\begin{align*}
    \{n\in\mathbb{N}:\lvert x_n-x_0\rvert<\epsilon\}\in\mathcal{F}.
\end{align*}
\end{definition}
Let us state a few well-known results pertaining to the convergence of a sequence under a filter without proof.
\begin{theorem}
{Let} $\mathcal{F}$ is a filter on $\mathbb{N}$ and $(x_n)_{n\in\mathbb{N}}$ and $(y_n)_{n\in\mathbb{N}}$ {be} two sequences in $\mathbb{K}.$    {Then the following hold true:}\\
$(i)$ $\lim\limits_{\mathcal{F}}x_n$, if exists, is unique.\\
$(ii)$ If $x=\lim\limits_{\mathcal{F}}x_n$, and $f:\mathbb{K}\to\mathbb{K}$ is continuous, {then} $f(x)=\lim\limits_{\mathcal{F}}f\left(x_n\right)$.\\
$(iii)$ If $x=\lim\limits_\mathcal{F}x_n$ and $y=\lim\limits_\mathcal{F}y_n$, then $x+\lambda y=\lim\limits_\mathcal{F}{(}x_n+\lambda y_n{)},$ for any $\lambda\in\mathbb{K}$. Also, $xy=\lim\limits_\mathcal{F}x_ny_n$ and $\frac{x}{y}=\lim\limits_\mathcal{F}\frac{x_n}{y_n}$ if $y_n\neq0\neq y,$ for every $n\in\mathbb{N}$.\\ 
$(iv)$ If $\mathcal{F}$ is an ultrafilter and $(x_n)_{n\in\mathbb{N}}$ is bounded, then $\lim\limits_\mathcal{F}x_n$ exists.
\end{theorem}
It is trivial to see that if $\mathcal{U}$ is the principal ultrafilter containing $\{N\},$ for some $N\in\mathbb{N}$, $\lim\limits_{\mathcal{U}}x_n=x_N$. We establish the {following} result pertaining to the limit of a bounded sequence under any free ultrafilter. 
\begin{prop}\label{subseq}
    {Let} $(x_n)_{n\in\mathbb{N}}$ and $(y_n)_{n\in\mathbb{N}}$ {be} two bounded sequences in $\mathbb{K}.$  {Then the following hold:} \\
    $(i)$ $x_0=\lim\limits_\mathcal{U}x_n$, for some free ultrafilter $\mathcal{U}$ on $\mathbb{N}$ if and only if $x_0$ is a subsequential limit of $x_n$.\\
    $(ii)$ $x_0=\lim\limits_\mathcal{U}x_n$ and $y_0=\lim\limits_\mathcal{U}y_n$, for some free ultrafilter $\mathcal{U}$ on $\mathbb{N}$ if and only if there exists an increasing sequence $(n_k)_{k\in\mathbb{N}}$ on $\mathbb{N}$ such that $x_{n_k}\to x_0$ and $y_{n_k}\to y_0$ as $k\to\infty$.
\end{prop}
\begin{proof}
Recall that an ultrafilter is free if and only if it contains no finite subset of $\mathbb{N}$.\\
$(i)$ We first prove the necessity. Suppose $x_0=\lim\limits_\mathcal{U}x_n,$ for some free ultrafilter $\mathcal{U}$ on $\mathbb{N}$. If $x_0$ is not a subsequential limit of $x_n$, there exists $\delta>0$ such that $\{n\in\mathbb{N}:\lvert x_n-x_0\rvert<\delta\}$ is finite {which is} a contradiction since $\mathcal{U}$ is a free ultrafilter.\par 
Now, assume {that} $(n_k)_{k\in\mathbb{N}}$ is an increasing sequence in $\mathbb{N}$ such that $x_{n_k}\to x_0,$ for some $x_0\in\mathbb{K}$ and consider the following set:
\begin{align}\label{base}
    \{A\subset\mathbb{N}:\{n_k:k\in\mathbb{N}\}\backslash A~\textit{is~finite}\}.
\end{align}
{Observe that the set defined by \eqref{base}} is a filter base and therefore is contained in some ultrafilter $\mathcal{U}$ on $\mathbb{N}$. Clearly, no finite subset of $\mathbb{N}$ is an element of $\mathcal{U}$ and hence $\mathcal{U}$ is an ultrafilter. Since, $\lim\limits_\mathcal{U}x_n=x_0$, the sufficiency is established.\\
$(ii)$ The sufficiency can be proved like  $(i)$. To prove the necessity, let there be no sequence $(n_k)_{k\in\mathbb{N}}$ in $\mathbb{N}$ such that $x_{n_k}\to x_0$ and $y_{n_k}\to y_0$. Then there exist $\delta_1,\delta_2>0$ such that
\begin{align*}
    \{n\in\mathbb{N}:\lvert x_n-x_0\rvert<\delta_1\}\cap\{n\in\mathbb{N}:\lvert y_n-y_0\rvert<\delta_2\}~\textit{is~finite},
\end{align*}
 a contradiction since $\mathcal{U}$ is a free ultrafilter and therefore contains no finite subset of $\mathbb{N}.$ {This proves} the necessity. 
\end{proof}

\subsection{Birkhoff-James orthogonality and smoothness of a point in $\ell_\infty$}\hfill
\\

We begin {this sub-section} by recalling a few known results:
\begin{theorem}\label{1}
The space $\ell_\infty$ is isometrically isomorphic to $C(\beta\mathbb{N})$, the Banach space of all $\mathbb{K}$-valued continuous functions on $\beta\mathbb{N}$ equipped with the supremum norm, where $\beta\mathbb{N}$ denotes the Stone-\v Cech compactification of $\mathbb{N}$. Recalling the homeomorphism between $\beta\mathbb{N}$ and the space of all ultrafilters on $\mathbb{N}$ equipped with the Stone topology, one can explicitly write down such an isometric isomorphism $T:\ell_\infty\to C(\beta\mathbb{N}),$ given by
\begin{align}\label{isomisom}
    T\left((x_n)_{n\in\mathbb{N}}\right)\left(\mathcal{U}\right)=\lim\limits_{\mathcal{U}}x_n,~~\mathcal{U}~\textit{an~ultrafilter~on~}\mathbb{N}.
\end{align}
\end{theorem}
Also, since $\beta\mathbb{N}$ is compact Hausdorff, by an application of the Riesz representation {Theorem} in measure theory, we {have} the following result:
\begin{theorem}\label{2}
    The dual space of $C(\beta\mathbb{N})$ is isometrically isomorphic to the space of all regular $\mathbb{K}$-valued Borel measures on $\beta\mathbb{N}$ equipped with the total variation norm and the functional corresponding to a regular $\mathbb{K}$-valued Borel measure $\mu$ acting on $C(\beta\mathbb{N})$ by
    \begin{align*}
        \mu:~f\mapsto\int\limits_{\beta\mathbb{N}}fd\mu,~~f\in C(\beta\mathbb{N}).
    \end{align*}
\end{theorem}
We note that by Theorem \ref{1} and {Theorem} \ref{2}, $\ell_\infty^*$ is isometrically isomorphic to the space of all regular $\mathbb{K}$-valued Borel measures on $\beta\mathbb{N}$ equipped with the total variation norm. We begin by characterizing the support functionals of a non-zero $f\in C(\beta\mathbb{N})$ and introduce the following definition in this regard.
\begin{definition}
For a given $f\in C(\beta\mathbb{N})$, we define $M_f$ to be the collection of all the points in $\beta\mathbb{N}$ where $f$ attains its norm, i.e.,
\begin{align*}
    M_f:=\{\mathcal{U}\in \beta\mathbb{N}:\lvert f(\mathcal{U})\rvert=\|f\|\}.
\end{align*}
\end{definition}
Using {the above}, we now characterize $J(f),$ (see \eqref{support} for definition) for a non-zero $f\in C(\beta\mathbb{N})$.
\begin{theorem}\label{support1}
    Let $f\in C(\beta\mathbb{N})$ be non-zero. Then $\mu\in J(f),$ for some regular $\mathbb{K}$-valued Borel measure $\mu$ if and only if 
    \begin{align*}
      \lvert \mu\rvert\left(\beta\mathbb{N}\setminus M_f\right)=0,~~\lvert \mu\rvert (M_f)=1~~\textit{and}~~d\mu(\mathcal{U})=\overline{\sgn(f(\mathcal{U}))}d\lvert \mu\rvert (\mathcal{U}),
    \end{align*}
    for almost every $\mathcal{U}\in M_f$ with respect to the measure $\mu$, where $\lvert \mu\rvert $ denotes the total variation of $\mu$.
\end{theorem}
\begin{proof}
The sufficiency follows by elementary computations. Now, if $\mu\in J(f)$, then $d\mu(\mathcal{U})=e^{i\theta(\mathcal{U})}d\lvert\mu\rvert(\mathcal{U}),$ for some {measurable function} $\theta:\beta\mathbb{N}\to\mathbb{R}$ . Note that
\begin{align*}
    \|f\|=\int\limits_{\beta\mathbb{N}}f(\mathcal{U})e^{i\theta(\mathcal{U})}d\lvert\mu\rvert(\mathcal{U})\leq\int\limits_{\beta\mathbb{N}}\lvert f(\mathcal{U})\rvert d\lvert\mu\rvert(\mathcal{U})\leq\|f\|.
\end{align*}
Hence equality must hold in both the inequalities involved. Equality in the second inequality {implies that} $\lvert f(\mathcal{U})\rvert=\|f\|,$ for almost every $\mathcal{U}\in\beta\mathbb{N},$ giving $\lvert \mu\rvert \left(\beta\mathbb{N}\setminus M_f\right)=0$ (and hence, $\lvert \mu\rvert (M_f)=1$) and  {equality in the} first inequality gives that for almost every $\mathcal{U}\in M_f,$ with respect to $\mu$, 
\begin{align*}
    f(\mathcal{U})e^{i\theta(\mathcal{U})}=\|f\|~~\Rightarrow~~e^{i\theta(\mathcal{U})}=\overline{\sgn(f(\mathcal{U}))}.
\end{align*}
\end{proof}

We now come to our characterization of Birkhoff-James orthogonality in $\ell_\infty$.
\begin{theorem}\label{orthogonality}
    {Let} $x=(x_n)_{n\in\mathbb{N}}$ and $y=(y_n)_{n\in\mathbb{N}}$ {be} two elements of $\ell_\infty.$  {Then}  the following are equivalent:\\
    $(i)$ $x\perp_By$.\\
    $(ii)$ $0\in\conv\left\{\lim\limits_\mathcal{U}\overline{x_n}y_n:\lim\limits_\mathcal{U}\lvert x_n\rvert=\|x\|\right\}$.\\
    $(iii)$ $0$ lies in the convex hull of $\{\overline{x_n}y_n:\lvert x_n\rvert =\|x\|\}$ and $\{\lim\limits_{k\to\infty}\overline{x_{n_k}}y_{n_k}:n_k~\textit{{is an }increasing~sequence~in~}\mathbb{N},~\lim\limits_{k\to\infty}\left\lvert x_{n_k}\right\rvert =\|x\|,~y_{n_k}~\textit{converge~as~}k\to\infty\}$
\end{theorem}
\begin{proof}
The result holds trivially if $x=0$. Hence, we assume that $x\neq0$.\\
$(i)\Leftrightarrow(ii)$\\
If $0\in\conv\{\lim\limits_\mathcal{U}\overline{x_n}y_n:\lim\limits_\mathcal{U}\lvert x_n\rvert=\|x\|\}$, then there exist  {ultrafilters $\mathcal{U}_i$} on $\mathbb{N}$ and $\lambda_i\in[0,1],$ for $1\leq i\leq m$ such that $\sum\limits_{i=1}^m\lambda_i=1$ and 
\begin{align*}
    \sum\limits_{i=1}^m\lambda_i\lim\limits_{\mathcal{U}_i}\overline{x_n}y_n=0,~~\lim\limits_{\mathcal{U}_i}\lvert x_n\rvert=\|x\|,~~\textit{for~every~}1\leq i\leq m.
\end{align*}
Consider the functional $\Psi:\ell_\infty\to\mathbb{K},$ given by\\
\begin{align*}
    \Psi\left((z_n)_{n\in\mathbb{N}}\right)=\sum\limits_{i=1}^n\lambda_i\lim\limits_{\mathcal{U}_i}\overline{\sgn(x_n)}z_n,~~(z_n)_{n\in\mathbb{N}}\in\ell_\infty.
\end{align*}
Then {clearly} $\Psi$ has norm 1 and $\Psi(x)=\|x\|$. Hence $\Psi$ is a support functional of $x$ that annihilates $y$ {and establishes}  the sufficiency.\par
Now, let us assume that $x\perp_By$. Then by Theorem \ref{1}, {we have} the two maps $f,g:\beta\mathbb{N}\to\mathbb{K},$ given by
\begin{align*}
    f(\mathcal{U}):=\lim\limits_{\mathcal{U}}x_n,~~g(\mathcal{U}):=\lim\limits_{\mathcal{U}}y_n,~~\mathcal{U}\in\beta\mathbb{N},
\end{align*}
satisfying $f\perp_Bg$ in $C(\beta\mathbb{N})$. Observe that $M_f=\{\mathcal{U}\in\beta\mathbb{N}:\lim\limits_\mathcal{U}\lvert x_n\rvert=\|x\|\}$. Now, by Theorem \ref{support1}, there must exist a regular positive Borel measure $\nu$ on $M_f$ with $\nu(M_f)=1$ such that
\begin{align*}
    \int\limits_{M_f}\overline{sgn(f(\mathcal{U}))}g(\mathcal{U})d\nu(\mathcal{U})=0~~\Leftrightarrow~~ \int\limits_{M_f}\overline{f(\mathcal{U})}g(\mathcal{U})d\nu(\mathcal{U})=0.
\end{align*}
Let $\Lambda$ be the collection of all regular positive Borel measures $\mu$ on $M_f$ with $\mu(M_f)=1$. Consider the map $\Phi:\Lambda\to\mathbb{K},$ given by
\begin{align*}
    \Phi(\mu):=\int\limits_{M_f}\overline{f(\mathcal{U})}g(\mathcal{U})d\mu(\mathcal{U}),~~ \mu\in\Lambda.
\end{align*}
Then clearly $\Phi(\Lambda)$ must be a convex subset of $\mathbb{K}$ since $\Lambda$ is convex. Observe that by Theorem \ref{support1}, $\Lambda$ is the collection of all support functionals of $\lvert f\rvert\in C(\beta\mathbb{N})$ and hence {is} weak* compact by the Banach-Alaoglu {Theorem}. Since the map $\Phi$ is continuous under the weak* topology, $\Phi(\Lambda)$ {will} be compact. Hence by the Krein-Milman {Theorem}, $\Phi(\Lambda)$ must be the closed convex hull of its extreme points.\par
We claim that  the only extreme points of the set $\Phi(\Lambda)$ are of the form $\overline{f(\mathcal{U})}g(\mathcal{U}),$ for some $\mathcal{U}\in M_f$. Let $\Phi(\mu)$ be an extreme point of $\Phi(\Lambda)$ where $\mu$ is not a Dirac delta measure at any point on $M_f$. Clearly, if $\overline{f}g$ is constant at every point of the support of $\mu$, then $\Phi(\mu)=\overline{f(\mathcal{U})}g(\mathcal{U}),$ for any $\mathcal{U}$ in the support of $\mu$. Therefore, let us assume that $\mathcal{U},\;\mathcal{V}$ are two points in the support of $\mu$ with $\overline{f(\mathcal{U})}g(\mathcal{U})\neq\overline{f(\mathcal{V})}g(\mathcal{V})$. Fix $0<\epsilon<\frac{1}{2}\lvert \overline{f(\mathcal{U})}g(\mathcal{U})-\overline{f(\mathcal{V})}g(\mathcal{V})\rvert$ and set:
\begin{align*}
    G_\epsilon:=\left\{\mathcal{W}\in M_f:\left\lvert \overline{f(\mathcal{U})}g(\mathcal{U})-\overline{f(\mathcal{W})}g(\mathcal{W})\right\rvert <\epsilon\right\}.
\end{align*}
Then $G_\epsilon$ is an open subset of $M_f$ containing $\mathcal{U}$ and $M_f\backslash G_\epsilon$ contains a neighbourhood of $\mathcal{V}$ in $M_f$. Hence $\mu(G_\epsilon),\mu(M_f\backslash G_\epsilon)>0$. Now since $\mu$ can be written as a convex combination of $\frac{1}{\mu\left(G_\epsilon\right)}\mu\vert_{G_\epsilon}$ and $\frac{1}{\mu\left(M_f\backslash G_\epsilon\right)}\mu\vert_{M_f\backslash G_\epsilon},$
\begin{align*}
    \Phi(\mu)=\frac{1}{\mu(G_\epsilon)}\int\limits_{G_\epsilon}\overline{f(\mathcal{W})}g(\mathcal{W})d\mu(\mathcal{W}),~~\textit{for~every~}0<\epsilon<\frac{1}{2}\lvert\overline{f(\mathcal{U})}g(\mathcal{U})-\overline{f(\mathcal{V})}g(\mathcal{V})\rvert,
\end{align*}
as $\Phi(\mu)$ is an extreme point of $\Phi(\Lambda)$. Hence,
\begin{align*}
    \left\lvert \Phi(\mu)-\overline{f(\mathcal{U})}g(\mathcal{U})\right\rvert &=\frac{1}{\mu(G_\epsilon)}\left\lvert \int\limits_{G_\epsilon}\left(\overline{f(\mathcal{W})}g(\mathcal{W})-\overline{f(\mathcal{U})}g(\mathcal{U})\right)d\mu(\mathcal{W})\right\rvert \\
    &\leq\frac{1}{\mu(G_\epsilon)}\int\limits_{G_\epsilon}\left\lvert \overline{f(\mathcal{W})}g(\mathcal{W})-\overline{f(\mathcal{U})}g(\mathcal{U})\right\rvert d\mu\left(\mathcal{{W}}\right)\leq\epsilon,
\end{align*}
giving $\Phi(\mu)=\overline{f(\mathcal{U})}g(\mathcal{U})$ since $0<\epsilon<\frac{1}{2}\lvert \overline{f(\mathcal{U})}g(\mathcal{U})-\overline{f(\mathcal{V})}g(\mathcal{V})\rvert$ is arbitrary.\par
Therefore, 
\begin{align*}
    0\in\overline{\conv}\left\{\lim\limits_\mathcal{U}\overline{x_n}y_n:\lim\limits_\mathcal{U}\lvert x_n\rvert=\|x\|\right\},{~~as~~ \Phi(\nu)=0}.
\end{align*}
Clearly, since $M_f$ is compact, $\{\lim\limits_\mathcal{U}\overline{x_n}y_n:\lim\limits_\mathcal{U}\lvert x_n\rvert=\|x\|\}$ is also compact. We now prove that the convex hull of a compact subset of $\mathbb{K}$ is closed. Indeed, using Caratheodory's theorem, for any element $x$ in $\conv(K),$ {where} $K\subset\mathbb{K}$ {is} compact, there exist $x_1,x_2,x_3\in K$ and $\lambda_1,\lambda_2,\lambda_3\in[0,1]$ such that
\begin{align*}
    \sum\limits_{i=1}^3\lambda_ix_i=x~~\textit{and}~~\sum\limits_{i=1}^3\lambda_i=1.
\end{align*}
Hence if $x_0\in\overline{\conv(K)}$, there exist sequences $x_i^{(n)}\in K$ and $\lambda_i^{(n)}\in[0,1],$ for $1\leq i\leq3$ and $n\in\mathbb{N}$ such that
\begin{align*}
    \lim\limits_{n\to\infty}\sum\limits_{i=1}^3\lambda_i^{(n)}x_i^{(n)}=x,~~\textit{and}~~\sum\limits_{i=1}^3\lambda_i^{(n)}=1,~~\textit{for~ every~}n\in\mathbb{N}.
\end{align*}
Since $K$ and $[0,1]$ are compact, considering a convergent subsequence of all the six sequences and passing on to the limits, {we have} $x_0\in\conv(K)$. Combining all the results, we therefore obtain
\begin{align*}
    0\in\conv\left\{\lim\limits_\mathcal{U}\overline{x_n}y_n:\lim\limits_\mathcal{U}\lvert x_n\rvert=\|x\|\right\}.
\end{align*}
$(ii)\Leftrightarrow(iii)$
This is immediate from Proposition \ref{subseq} and the fact that if $\mathcal{U}$ is a principal ultrafilter containing $\{N\},$ for some $N\in\mathbb{N}$, then $\lim\limits_\mathcal{U}x_n=x_N,$ for every sequence $(x_n)_{n\in\mathbb{N}}$ in $\ell_\infty$.
\end{proof}
We record the characterization of Birkhoff-James orthogonality in $\ell_\infty,$ for the real case as a corollary below.
\begin{cor}
    {Let} $x=(x_n)_{n\in\mathbb{N}}$ and $y=(y_n)_{n\in\mathbb{N}}$ be two elements of $\ell_\infty$ over $\mathbb{R}$. {Then} $x\perp_By$ if and only if any of the {following holds true} :\\
    $(i)$ There is an increasing sequence of natural numbers $(n_k)_{k\in\mathbb{N}}$ such that $\left\lvert x_{n_k}\right\rvert \to\|x\|$ and $x_{n_k}y_{n_k}\to0$ as $k\to\infty$. \\
    $(ii)$ There are increasing sequences of natural numbers $(n_k)_{k\in\mathbb{N}}$ and $(m_k)_{k\i\mathbb{N}}$ such that $\left\lvert x_{n_k}\right\rvert \to\|x\|$ and $\left\lvert x_{m_k}\right\rvert \to\|x\|$ as $k\to\infty$ with $x_{n_k}y_{n_k}\geq0\geq x_{m_k}y_{m_k},$ for every $k\in\mathbb{N}$.\\
    $(iii)$ There is an increasing sequence of natural numbers $(n_k)_{k\in\mathbb{N}}$ such that $\left\lvert x_{n_k}\right\rvert \to\|x\|$ as $k\to\infty$ and a natural number $N$ such that $\lvert x_N\rvert=\|x\|$ and $x_Ny_N$ and $x_{n_k}y_{n_k}$ are of different signs for every $k\in\mathbb{N}$.\\
    $(iv)$ There are $N,M\in\mathbb{N}$ such that $\lvert x_N\rvert=\|x\|=\lvert x_M\rvert$ and $x_Ny_N\geq0\geq x_My_M$.
\end{cor}
We conclude this section by characterizing the smooth points of $\ell_\infty$.
\begin{theorem}\label{smooth}
    {Let} $x=(x_n)_{n\in\mathbb{N}}$ be a non-zero element in $\ell_\infty$. {Then} $x$ is a smooth point if and only if there is no subsequence of $(\lvert x_n\rvert)_{n\in\mathbb{N}}$ that converges to $\|x\|$ and there exists unique $N\in\mathbb{N}$ such that $\lvert x_N\rvert=\|x\|$.
\end{theorem}
\begin{proof}
{We first prove the sufficiency. Suppose, there is no subsequence of $(\lvert x_n\rvert)_{n\in\mathbb{N}}$ that converges to $\|x\|$ and there exists a unique $N\in\mathbb{N}$ such that $\lvert x_N\rvert=\|x\|$. Then it follows from Theorem \ref{orthogonality} that $x\perp_By,$ for some $y=(y_n)_{n\in\mathbb{N}}\in\ell_\infty$, if and only if $y_N=0.$}  {Consequently,} $x\perp_By$ and $x\perp_Bz$ implies $x\perp_B(y+z),$ proving $x$ is smooth.\par
Again, if $x\in\ell_\infty$ is smooth, then the function $f:\beta\mathbb{N}\to\mathbb{K},$ given by
\begin{align*}
    f(\mathcal{U}):=\lim\limits_\mathcal{U}x_n,~~\mathcal{U}\in\beta\mathbb{N},
\end{align*}
is smooth in $C(\beta\mathbb{N})$ (Theorem \ref{1}). Therefore, if $\mathcal{U},\mathcal{V}\in M_f$, then  $\Psi,\Phi:C(\beta\mathbb{N})\to\mathbb{K},$ given by
\begin{align*}
    \Psi(g):=\overline{\sgn(f(\mathcal{U}))}g(\mathcal{U}),~~\Phi(g):=\overline{\sgn(f(\mathcal{V}))}g(\mathcal{V}),~~g\in C(\beta\mathbb{N}),
\end{align*}
are two distinct support functionals of $f$ because there exists a continuous map $F:\beta\mathbb{N}\to[0,1]$ such that $F(\mathcal{U})=0$ and $F(\mathcal{V})=1,$ since $\beta\mathbb{N}$ is compact Hausdorff. Hence $M_f$ must be singleton.\par
Now, if $\mathcal{U}\in M_f$ is a free ultrafilter, then by Theorem \ref{subseq}, there exists an increasing sequence $(n_k)_{k\in\mathbb{N}}$ of natural numbers such that $\left\lvert x_{n_k}\right\rvert \to \|x\|$ as $k\to\infty$. Consider the following two collections:
\begin{align*}
    \{A\subset \mathbb{N}:\{n_{2k}:k\in\mathbb{N}\}\backslash A~\textit{is~finite~}\}~~\textit{and}~~\{A\subset \mathbb{N}:\{n_{2k-1}:k\in\mathbb{N}\}\backslash A~\textit{is~finite~}\}.
\end{align*}
Clearly, both {of these} sets are filter bases on $\mathbb{N}$ and are therefore contained in free ultrafilters. Also note that since the first collection contains $\{n_{2k}:k\in\mathbb{N}\}$ and the second collection contains $\mathbb{N}\backslash\{n_{2k}:k\in\mathbb{N}\}$, the two ultrafilters must be distinct. However, the sequence {(}$\left\lvert x_n\right\rvert ${)} has the same limit under both {of} the ultrafilters. Hence if $M_f$ is singleton, it can contain only a principle ultrafilter and {thus} the necessity {follows}.
\end{proof}

\subsection{Local symmetry of Birkhoff-James orthogonality in $\ell_\infty$}\hfill
\\

In this sub-section we characterize the left-symmetric and the right-symmetric points of the space $\ell_\infty$.
\begin{theorem}\label{ls}
    The only non-zero left-symmetric points of $\ell_\infty$ are scalar multiples of $e_n$ for $n\in\mathbb{N}$, where $e_n$ denotes the sequence having the $n$-th term 1 and the rest of the terms 0.
\end{theorem}
\begin{proof}
Since $\lim\limits_\mathcal{U}e_n=0,$ for every ultrafilter other than the principal ultrafilter containing $\{n\}$, by Theorem \ref{orthogonality}, $e_n\perp_Bx,$ for some $x=(x_n)_{n\in\mathbb{N}}\in\ell_\infty$ if and only if $x_n=0$. Hence for $\lambda\in\mathbb{K}$,
\begin{equation*}
    \|x+\lambda e_n\|=\max\left\{\lvert\lambda\rvert,\|x\|\right\}\geq\|x\|,
\end{equation*}
{and this proves} the necessity.\par
Now, let $x=(x_n)_{n\in\mathbb{N}}$ be a left-symmetric point of $\ell_\infty$. If $\lim\limits_{\mathcal{U}}\lvert x_n\rvert=\|x\|,$ for some free ultrafilter $\mathcal{U}$ on $\mathbb{N}$, then set  $y=(y_n)_{n\in\mathbb{N}}\in\ell_\infty,$ given by $y_n:=\frac{1}{n}\sgn(x_n)$. Then clearly, $\lim\limits_\mathcal{U}\overline{x_n}y_n=0,$ since $\mathcal{U}$ is a free ultrafilter. Hence, $x\perp_By.$  However, by Theorem \ref{orthogonality}, for some $z=(z_n)_{n\in\mathbb{N}}\in\ell_\infty$, $y\perp_Bz$ if and only if $z_N=0$ where $N=\min\{n\in\mathbb{N}:x_n\neq0\}$. {Thus} $y\not\perp_Bx,$ {implying} that $\lim\limits_\mathcal{U}\lvert x_n\rvert\neq\|x\|,$ for any free ultrafilter $\mathcal{U}$ on $\mathbb{N}$.\par
Hence, there exists $N\in\mathbb{N}$ such that $\lvert x_N\rvert =\|x\|$. Now, suppose {that} $x_M\neq0,$ for some $M\neq N\in \mathbb{N}$. Then clearly, setting $y=e_M$, yields that $x\perp_By$ and $y\not\perp_Bx,$ by Theorem \ref{orthogonality}, {which establishes}  the necessity.
\end{proof}
Now, from Corollary \ref{isometry}, if $T:\ell_\infty\to\ell_\infty$ is an onto linear isometry, then for every $n\in\mathbb{N}$, $T(e_n)=ce_m,$ for some $m\in\mathbb{N}$ and {some unimodular constant $c$.} Also, since $T$ is onto, $T$ is invertible and $T^{-1}$ is also an onto isometry and hence $T$ must be a signed permutation operator. Note that this extends the result of Lamperti \cite{BL} on  $L_p$ spaces for $1\leq p<\infty$ to the $p=\infty$ case for onto operators with the measure space as $\mathbb{N}$ under counting measure. We record this result as a corollary.
\begin{cor}\label{isom}
{Let} $T:\ell_\infty\to\ell_\infty$ {be} an onto isometry.  {Then}  $T$ must be a signed permutation operator.
\end{cor}
We conclude this sub-section by characterizing the right-symmetric points of $\ell_\infty$.
\begin{theorem}\label{rs}
    $x=(x_n)_{n\in\mathbb{N}}\in\ell_\infty$ is a right-symmetric point if and only if $\lvert x_n\rvert=\|x\|$ for every $n\in\mathbb{N}$.
\end{theorem}
\begin{proof}
Note that if $\lvert x_n\rvert=\|x\|,$ for every $n\in\mathbb{N}$, then $\lim\limits_\mathcal{U}\lvert x_n\rvert=\|x\|,$ for every ultrafilter $\mathcal{U}$ on $\mathbb{N}$ and hence by Theorem \ref{orthogonality}, $y\perp_B x,$ for $y=(y_n)_{n\in\mathbb{N}}\in\ell_\infty$ if and only if
\begin{align*}
    & 0\in \conv\left\{\lim\limits_\mathcal{U}\overline{y_n}x_n:\lim\limits_\mathcal{U}\lvert y_n\rvert=\|y\|\right\}\\
    \Rightarrow~~ & 0\in \conv\left\{\lim\limits_\mathcal{U}\overline{x_n}y_n:\lim\limits_\mathcal{U}\lvert x_n\rvert=\|x\|\right\}\\
    \Rightarrow ~~& x\perp_B y.
\end{align*}
{So we establish} the necessity.\par
Now, if $x\neq0$ is right-symmetric and $\lvert x_N\rvert<\|x\|$ for some $N\in\mathbb{N}$, set $y=(y_n)_{n\in\mathbb{N}}\in\ell_\infty$ given by 
\begin{align*}
    y_n:=
    \begin{cases}
    \sgn(x_n),~~n\neq N,\\
    -\sgn(x_N),~~n=N~\textit{and~if~}x_N\neq0,\\
    1,~~n=N~\textit{and~if~}x_N=0.
    \end{cases}
\end{align*}
Then clearly, by Theorem \ref{orthogonality}, $y\perp_Bx$. However,
\begin{align*}
    \overline{x_n} y_n:=
    \begin{cases}
  \lvert x_n\rvert ,~~n\neq N,\\
    -\lvert x_N\rvert ,~~n=N.
    \end{cases}
\end{align*}
Hence, if $\mathcal{U}$ is an ultrafilter not containing $\{N\}$,
\begin{align}\label{limit}
    \lim\limits_\mathcal{U}\overline{x_n}y_n=\lim\limits_\mathcal{U}\lvert x_n\rvert.
\end{align}
Therefore, since $\lvert x_N\rvert<\|x\|$, the limit of $\lvert x_n\rvert$ under the principal ultrafilter containing $\{N\}$ is not $\|x\|$ and hence by \eqref{limit} and Theorem \ref{orthogonality}, $x\not\perp_By$ establishing the necessity.

\end{proof}

\subsection{Geometry of $c$, $c_0$ and $c_{00}$}\hfill
\\

Recall that $c$, $c_0$ and $c_{00}$ are the collections of all convergent sequences, convergent to zero sequences and eventually zero sequences respectively. Let us denote the limit of a sequence $x=(x_n)_{n\in\mathbb{N}}\in c$ by $\lim x$.
Since all the three spaces are subspaces of $\ell_\infty$, two elements in any of these spaces are Birkhoff-James orthogonal if and only if they are Birkhoff-James orthogonal in $\ell_\infty$. Keeping this fundamental principle in mind {we observe} that {if} $x=(x_n)_{n\in\mathbb{N}}\in c$, $\lim\limits_\mathcal{U}x_n=\lim x$ and {if} $x\in c_0$, $\lim\limits_\mathcal{U}x_n=0,$ for every free ultrafilter $\mathcal{U}$ on $\mathbb{N}$. {Thus} we obtain the following result from Theorem \ref{orthogonality}.
\begin{theorem}\label{Orthogonality subcases}
Suppose {that} $x=(x_n)_{n\in\mathbb{N}}$ and $y=(y_n)_{n\in\mathbb{N}}$ are two sequences in $\mathbb{K}$. Then\\
    $(i)$ If $x,y\in c$, then $x\perp_By$ if and only if 
    \begin{align*}
        0\in&\conv\left(\{\overline{x_n}y_n:\lvert x_n\rvert=\|x\|\}\cup\{\lim \overline{x}y\}\right),~~\textit{if~}\lim\lvert x\rvert=\|x\|,\\
        \textit{and~}~0\in&\conv\{\overline{x_n}y_n:\lvert x_n\rvert=\|x\|\},~~\textit{if~}\lim\lvert x\rvert\neq\|x\|.
    \end{align*}
    In particular, for $\mathbb{K}=\mathbb{R}$, $x\perp_By$ if and only if any of the following is true:\\
    \begin{enumerate}
        \item There exists $N\in\mathbb{N}$ such that $\lvert x_N\rvert=\|x\|$ and $y_N=0$ or $\lim\lvert x\rvert=\|x\|$ and $\lim y=0$.
        \item There exists $N\in\mathbb{N}$ such that $\lvert x_N\rvert=\lim\lvert x\rvert=\|x\|$ and $x_Ny_N$ and $\lim xy$ are of different signs.
        \item There exist $N,M\in\mathbb{N}$ such that $\lvert x_N\rvert=\lvert x_M\rvert=\|x\|$ and $x_Ny_N<0<x_My_M$.
    \end{enumerate}
    $(ii)$ If $x,y\in c_0$ or $c_{00}$, $x\perp_By$ if and only if
    \begin{align*}
        0\in&\conv\{\overline{x_n}y_n:\lvert x_n\rvert=\|x\|\}.
    \end{align*}
    In particular, for $\mathbb{K}=\mathbb{R}$, $x\perp_By$ if and only if one of the following holds:\\
    \begin{enumerate}
        \item There exists $N\in\mathbb{N}$ such that $\lvert x_N\rvert=\|x\|$ and $y_N=0$.
        \item There exist $N,m\in\mathbb{N}$ such that $\lvert x_N\rvert=\lvert x_M\rvert=\|x\|$ and $x_Ny_N<0<x_My_M$.
    \end{enumerate}
\end{theorem}
The characterization of the smooth points of these spaces requires a little bit more work.
\begin{theorem}\label{smooth subspace}
$(i)$ $x=(x_n)_{n\in\mathbb{N}}\in c$ is smooth if and only if there exists a unique $N\in\mathbb{N}$, such that $\lvert x_N\rvert=\|x\|$ and $\lim\lvert x\rvert\neq\|x\|$ or $\lvert x_n\rvert<\|x\|,$ for every $n\in\mathbb{N}$.\\
$(ii)$ $x=(x_n)_{n\in\mathbb{N}}\in c_0$ or $c_{00}$ is smooth if and only if there exists a unique $N\in\mathbb{N}$ such that $\lvert x_N\rvert=\|x\|$.
\end{theorem}
\begin{proof}
$(i)$ If there are distinct $N,M\in\mathbb{N}$ such that $\lvert x_N\rvert=\lvert x_M\rvert=\|x\|$, then the two functionals $\Psi,\Phi:c\to\mathbb{K}$ given by
\begin{align*}
    \Psi(y):=\overline{\sgn(x_N)}y_N,~~\Phi(y):=\overline{\sgn(x_M)}y_M,~~y=(y_n)_{n\in\mathbb{N}}\in c,
\end{align*}
are distinct support functionals of $x$. Again if $\lvert x_N\rvert=\lim\lvert x\rvert=\|x\|,$ for some $N\in\mathbb{N}$, then $\Psi',\Phi':c\to\mathbb{K}$ given by
\begin{align*}
    \Psi'(y):=\overline{\sgn(x_N)}y_N,~~\Phi'(y):=\lim\limits_{n\to\infty}\overline{\sgn(x_n)}y_n,~~y=(y_n)_{n\in\mathbb{N}}\in c,
\end{align*}
are two distinct support functionals of $x$, establishing the necessity.\par
Now, by Theorem \ref{Orthogonality subcases}, if there exists {a} unique $N\in\mathbb{N}$ with $\lvert x_N\rvert=\|x\|$ and $\lim\lvert x\rvert<\|x\|$, then $x\perp_By,$ for some $y=(y_n)_{n\in\mathbb{N}}\in c$ if and only if $y_N=0$ and {clearly}, Birkhoff-James orthogonality is right additive at $x$. Again if $\lvert x_n\rvert<\|x\|,$ for every $n\in\mathbb{N}$, then $\lim\lvert x\rvert=\|x\|$ and again by Theorem \ref{Orthogonality subcases}, $x\perp_B y,$ for some $y\in c$ if and only if $\lim y=0$ proving the right additivity of Birkhoff-James orthogonality at $x$ and  the sufficiency {is established}.\par
$(ii)$ The result for $c_0$ and $c_{00}$ follows from part $(i)$ of this {theorem} and the observation that for $x\in c_0$ or $c_{00}$, $\lim\lvert x\rvert=0\neq\|x\|$ unless $x=0$.
\end{proof}
Using these two results and going through the proofs of Theorem \ref{ls}, {Theorem} \ref{rs} and Corollary \ref{isom}, we obtain the following result.
\begin{theorem}
\begin{enumerate}
    \item The left-symmetric points of each of these spaces are $ce_n$ for $n\in\mathbb{N}$ and $c\in\mathbb{K}$.
    \item The right-symmetric points of $c$ are the sequences $x=(x_n)_{n\in\mathbb{N}}$ such that $\lvert x_n\rvert=\|x\|,$ for every $n\in\mathbb{N}$ and the spaces $c_0$ and $c_{00}$ have no non-zero right-symmetric point.
    \item The isometries of each of these spaces onto itself are the signed permutations operators. 
\end{enumerate}
     
\end{theorem}
\section{Geometry of $\ell_1$}
In this section, we characterize Birkhoff-James orthogonality and its local symmetry in $\ell_1$. We begin with a known result.
\begin{theorem}
The dual space of $\ell_1$ is isometrically isomorphic to $\ell_\infty$ with the functional $\Psi_b$ corresponding to an element $b=(b_n)_{n\in\mathbb{N}}\in\ell_\infty ,$ given by
\begin{align*}
    \Psi_b(a)=\sum\limits_{n=1}^\infty b_na_n,~~a=(a_n)_{n\in\mathbb{N}}\in\ell_1.
\end{align*}
\end{theorem}
In order to avoid confusion, we {denote} the norm of an element in $\ell_1$ by $\|.\|_1$ and the norm of an element in its dual, identified with the $\ell_\infty$ space by $\|.\|_\infty$. We come to the following preliminary lemma before giving a complete characterization of Birkhoff-James orthogonality in $\ell_1$.
\begin{lemma}\label{norming}
{Let $a =(a_n)_{n\in\mathbb{N}}\in\ell_1$ be non-zero. Then $b = (b_n)_{n\in \mathbb{N}}\in\ell_\infty$ is a support functional of  $a$} if and only if $b_n=\overline{\sgn(a_n)}$ if $a_n\neq0$ and $\lvert b_n\rvert\leq1$ if $a_n=0$.

\end{lemma}
\begin{proof}
The sufficiency can be established by elementary calculations. For the necessity, note that 
\begin{align*}
    \|a\|_1=\sum\limits_{n=1}^\infty b_na_n=\sum\limits_{a_n\neq0}b_na_n\leq\sum\limits_{a_n\neq0}\lvert b_na_n\rvert\leq\|b\|_\infty\sum\limits_{a_n\neq0}\lvert a_n\rvert=\|a\|_1.
\end{align*}
Hence, equality holds in both the inequalities, giving $\lvert b_n\rvert=\|b\|_\infty=1$ and $b_na_n=\lvert b_na_n\rvert$, i.e., $b_n=\overline{\sgn(a_n)},$ for every $n\in\mathbb{N}$ with $a_n\neq0$. Also, since $\|b\|_\infty=1$, $\lvert b_n\rvert\leq1,$ for $n\in\mathbb{N}$ with $a_n=0$.
\end{proof}
We now come to the characterization of Birkhoff-James orthogonality in $\ell_1$.
\begin{theorem}\label{orthogonality2}
{Let} $a=(a_n)_{n\in\mathbb{N}}$  {and} $b=(b_n)_{n\in\mathbb{N}}\in\ell_1.$  {Then}  $a\perp_Bb$ if and only if 
\begin{align}\label{orthogonal}
    \left\lvert \sum\limits_{n=1}^\infty\overline{\sgn(a_n)}b_n\right\rvert \leq\sum\limits_{a_n=0}\lvert b_n\rvert,
\end{align}
where sum over an empty set is defined to be 0.
\end{theorem}
\begin{proof}
{If} $a\perp_Bb$, then there exists $c=(c_n)_{n\in\mathbb{N}}\in\ell_\infty$ with $\|c\|_\infty=1$ such that $\Psi_c\in J(a)$ and $\Psi_c(b)=0$. This yields
\begin{align*}
    \left\lvert \sum\limits_{a_n\neq0}\overline{\sgn(a_n)}b_n\right\rvert \leq\sum\limits_{a_n=0}\lvert c_nb_n\rvert\leq\sum\limits_{a_n=0}\lvert b_n\rvert,
\end{align*}
proving the necessity.\par
Again, if \eqref{orthogonal} holds, set $k\in\mathbb{K}$, $\lvert k\rvert\leq1$ such that 
\begin{align*}
    \sum\limits_{n=1}^\infty\overline{\sgn(a_n)}b_n=k\left(\sum\limits_{a_n=0}\lvert b_n\rvert\right).
\end{align*}
Consider $c=(c_n)_{n\in\mathbb{N}}$ given by
\begin{align*}
    c_n:=
    \begin{cases}
    \overline{\sgn(a_n)},~~a_n\neq0,\\
    -k\,\overline{\sgn(b_n)},~~a_n=0.
    \end{cases}
\end{align*}
Then clearly $c\in\ell_\infty$ and by Lemma \ref{norming}, $\Psi_c\in J(a)$. Since $\Psi_c(b)=0$, the sufficiency follows.
\end{proof}
As a corollary to {the} above result, we mention the case {when} $\mathbb{K}=\mathbb{R}.$
\begin{cor}
Suppose {that $a=(a_n)_{n\in\mathbb{N}}$ {and} $b=(b_n)_{n\in\mathbb{N}}\in\ell_1$ are two members of $\ell_1$ over $\mathbb{R}$.} Let $\mathcal{N}_1:=\{n\in\mathbb{N}:a_nb_n>0\}$, $\mathcal{N}_2:=\{n\in\mathbb{N}:a_nb_n<0\}$ and $\mathcal{N}_0:=\{n\in\mathbb{N}:a_n=0\}$. Then $a\perp_Bb$ if and only if
\begin{align*}
    \left\lvert \sum\limits_{n\in\mathcal{N}_1}\lvert b_n\rvert-\sum\limits_{n\in\mathcal{N}_2}\lvert b_n\rvert\right\rvert \leq\sum\limits_{n\in\mathcal{N}_0}\lvert b_n\rvert.
\end{align*}
\end{cor}
We now come to the characterization of the smooth points of $\ell_1$.
\begin{theorem}
The smooth points of $\ell_1$ are the sequences having no zero term.
\end{theorem}
\begin{proof}
Suppose $a=(a_n)_{n\in\mathbb{N}}\in\ell_1$ has no zero term. Then by Theorem \ref{orthogonality2}, if $a\perp_Bb$ {and} {$a\perp_Bc$} where $b=(b_n)_{n\in\mathbb{N}},c=(c_n)_{n\in\mathbb{N}}\in\ell_1$, then 
\begin{align*}
    \sum\limits_{n=1}^\infty\overline{\sgn(a_n)}b_n=\sum\limits_{n=1}^\infty\overline{\sgn(a_n)}c_n=0.
\end{align*}
Hence
\begin{align*}
    \sum\limits_{n=1}^\infty\overline{\sgn(a_n)}(b_n+c_n)=0~~\Rightarrow~~a\perp_B(b+c),
\end{align*}
 proving $a$ is smooth.\par
 If $a=(a_n)_{n\in\mathbb{N}}\in\ell_1$ has a zero term $a_N,$ for some $N\in\mathbb{N}$, by Lemma \ref{norming}, {we obtain} $b=(b_n)_{n\in\mathbb{N}}$ {and} $c=(c_n)_{n\in\mathbb{N}}\in\ell_\infty,$ given by
\begin{align*}
    b_n:=
    \begin{cases}
    \overline{\sgn(a_n)},~~a_n\neq0,\\
    0,~~a_n=0,
    \end{cases}
    ~~\textit{and~~} c_n:=
    \begin{cases}
    \overline{\sgn(a_n)},~~a_n\neq0,\\
    0,~~a_n=0~\textit{and~}n\neq N,\\
    1,~~n=N.
    \end{cases}
\end{align*}
are two distinct support functionals of $a.$ Hence $a$ cannot be smooth.
\end{proof}
We finally come to the characterization of local symmetry of Birkhoff-James orthogonality in $\ell_1$.
\begin{theorem}
No non-zero point of $\ell_1$ is left-symmetric.
\end{theorem}
\begin{proof}
Suppose {that} $a=(a_n)_{n\in\mathbb{N}}\in\ell_1$ is non-zero. If there exists $N\in\mathbb{N}$ such that $a_N=0$, {then} setting $b=(b_n)_{n\in\mathbb{N}},$ given by $b_n:=\sgn(a_n)\frac{1}{2^n}$,  {we have} $\|b\|_1<\infty$. {Thus}  considering $c:=b+2\|b\|_1e_N,$  Theorem \ref{orthogonality2}  {implies} $a\perp_Bc$ and $c\not\perp_Ba$. If $a_n\neq0,$ for every $n\in\mathbb{N}$, then there exists $M\in\mathbb{N}$ such that
\begin{align*}
    \sum\limits_{n=1}^M\lvert a_n\rvert\neq\sum\limits_{n=M+1}^\infty\lvert a_n\rvert.
\end{align*}
Set $b=(b_n)_{n\in\mathbb{N}},$ given by,
\begin{align*}
    b_n:=
    \begin{cases}
    \sgn(a_n),~~1\leq n\leq M,\\
    -\sgn(a_n)\frac{M}{2^{n-M}},~~n\geq M+1.
    \end{cases}
\end{align*}
Then $b\in\ell_1.$  {Again,} by Theorem \ref{orthogonality2}, $a\perp_Bb$ and $b\not\perp_Ba$ and hence $a$ is not a left-symmetric point.
\end{proof}
We now come to the characterization of the right-symmetric points of $\ell_1$ and thereby find the collection of onto isometries of the space.
\begin{theorem}
    The only right-symmetric points of $\ell_1$ are scalar multiples of $e_n,$ for $n\in\mathbb{N}$.
\end{theorem}
\begin{proof}
Observe that for $b=\left(b_n\right)_{n\in\mathbb{N}}\in\ell_1$, by Theorem  \ref{orthogonality2}, $b\perp_B ce_n,$ for some $c\neq 0$ if and only if $b_n=0$ and hence $ce_n\perp_B b$. Again if $a=\left(a_n\right)_{\mathbb{N}}\in \ell_1$ has $a_n,a_m\neq0,$ for $n\neq m$, we consider $r\in\mathbb{N}$ such that 
\begin{equation*}
    0<\lvert a_r\rvert\leq\sum\limits_{k\neq r}\lvert a_k\rvert.
\end{equation*}
Hence by Theorem \ref{orthogonality2}, $a\not\perp_Be_r$ but $e_r\perp_B a$.
\end{proof}
We use this result to characterize the onto isometries of $\ell_1$ with the help of Corollary \ref{isometry} and thereby establish the Banach-Lamperti theorem for onto operators on $\ell_1$.
\begin{cor}
{Let} $T:\ell_1\to\ell_1$ {be} an onto isometry. {Then} $T$ must be a signed permutation operator.
\end{cor}

\section{Geometry of $\ell_p,$ for $1<p<\infty$, $p\neq2$}
In this final section, we characterize Birkhoff-James orthogonality and its local symmetry in $\ell_p,$ for $p\in(1,\infty)\backslash\{2\}$. We begin with  {an} well-known result.
\begin{theorem}
The dual space of $\ell_p$ is isometrically isomorphic to $\ell_q$ where $\frac{1}{p}+\frac{1}{q}=1$ with the functional $\Psi_a$ corresponding to $a=(a_n)_{n\in\mathbb{N}}\in\ell_q,$ given by
\begin{align*}
    \Psi_a(x)=\sum\limits_{n=1}^\infty a_nx_n,~~x=(x_n)_{n\in\mathbb{N}}\in\ell_p.
\end{align*}
\end{theorem}
Hence, as was done in the case of $\ell_1$, we {continue} to denote the norm of an element of $\ell_p$ by $\|.\|_p$ and the norm of an element in the dual of $\ell_p$, identified with $\ell_q$, by $\|.\|_q$.  {At the beginning, we note} down a corollary of this theorem pertaining to the characterization of the support functional of an element of $\ell_p$.
\begin{cor}\label{normin}
{Let} $a=(a_n)_{n\in\mathbb{N}}\in\ell_q$ {be} a support functional of $x=(x_n)_{n\in\mathbb{N}}\in\ell_p\setminus\{0\}$. {Then} $a_n=\frac{1}{\|x\|_p^{p-1}}\overline{\sgn(x_n)}\lvert x_n\rvert^{p-1}$, $n\in\mathbb{N}$.
\end{cor}
The proof of this result involves elementary computations and the  {equality criteria of}  Holder's inequality. \par
We now come to our characterization of Birkhoff-James orthogonality in $\ell_p$ which follows as a direct consequence of Corollary \ref{normin} and James' characterization of Birkhoff-James orthogonality.
\begin{theorem}
{Let} $x=(x_n)_{n\in\mathbb{N}}$  {and}     $y=(y_n)_{n\in\mathbb{N}}\in\ell_p.$  {Then}  $x\perp_By$ if and only if 
\begin{align*}
    \sum\limits_{n=1}^\infty \overline{\sgn(x_n)}\lvert x_n\rvert^{p-1}y_n=0.
\end{align*}
\end{theorem}
We now characterize the local symmetry of Birkhoff-James orthogonality in $\ell_p.$
\begin{theorem}\label{orthogonality3}
$x=(x_n)_{n\in\mathbb{N}}\in\ell_p$ is a left-symmetric point if and only if $x$ is a non-zero right-symmetric point if and only if $x=ce_N,$ for some $N\in\mathbb{N}$, $c\in\mathbb{K}$ or $x=c_1e_N+c_2e_M,$ for some $N,M\in\mathbb{N}$ and $c_1,c_2\in\mathbb{K}$, $\lvert c_1\rvert=\lvert c_2\rvert$.
\end{theorem}
\begin{proof}
We first characterize the left-symmetric points. Let $x\neq0$. The sufficiency can be verified by elementary computations. Now, suppose $x_N,x_M\neq0$ and $\lvert x_N\rvert\neq\lvert x_M\rvert,$ for some $N,M\in\mathbb{N}$. Set $y=\sgn(x_N)\lvert x_M\rvert^{p-1}e_N-\sgn(x_M)\lvert x_N\rvert^{p-1}$ and note that by Theorem \ref{orthogonality3}, $x\perp_By$. However,
\begin{align*}
    \sum\limits_{n=1}^\infty \overline{\sgn(y_n)}\lvert y_n\rvert^{p-1}x_n=\lvert x_N\rvert^{p^2-2p}-\lvert x_M\rvert^{p^2-2p}\neq0,
\end{align*}
since $p\neq2$, proving $y\not\perp_Bx$. Again if $x_N,x_M,x_K\neq0$, we may assume that $\lvert x_N\rvert=\lvert x_M\rvert=\lvert x_K\rvert,$ as otherwise, by the previous argument, $x$ cannot be left-symmetric. But then, setting $y=\sgn(x_N)e_N-\frac{1}{2}(\sgn(x_M)e_M-\sgn(x_K)e_K)$ clearly yields $x\perp_By$ and $y\not\perp_Bx$ by Theorem \ref{orthogonality3}.\par 
Now, Proposition 2.1 of \cite{KP} states that in a smooth, strictly convex space, a point is left-symmetric if and only if it is right-symmetric. Hence the necessity of the right-symmetric case follows from the left-symmetric case.
\end{proof}
As a consequence of this result, we characterize all the onto isometries $T:\ell_p\to\ell_p$ and prove the Banach-Lamperti {Theorem} for onto operators on $\ell_p$.
\begin{theorem}
{Let} $T:\ell_p\to\ell_p$ {be} an onto isometry. {Then} $T$ must be a signed permutation operator.
\end{theorem}
\begin{proof}
Observe that $T$ and $T^{-1}$ are both onto isometries and hence $T(x)$ and $T^{-1
}(x)$ are both left-symmetric points of $\ell_p,$ for any $x\in\ell_p,$  by Corollary \ref{isometry}. Hence it is sufficient to show that for every $n\in\mathbb{N}$, there exists $m\in\mathbb{N}$ such that $T(e_n)=ce_m,$ for some $c\in\mathbb{K}$ with $\lvert c\rvert =1$. Suppose by contradiction, $T(e_n)\neq ce_m,$ for any $m\in\mathbb{N}$ and $\lvert c\rvert =1$. Then $T(e_n)=\frac{1}{2^{\frac{1}{p}}}(c_1e_i+c_2e_j),$ for $i\neq j\in\mathbb{N}$ and $c-1,c_2\in\mathbb{K}$, $\lvert c_1\rvert=\lvert c_2\rvert=1$. Now, $T^{-1}\left(\frac{1}{2^{\frac{1}{p}}}(c_1e_i-c_2e_j)\right)$ is a left-symmetric point other than any unimodular multiple of $e_n$ and hence must be $ce_m$ for some $m\neq n$ and $\lvert c\rvert=1$ or $\frac{1}{2^{\frac{1}{p}}}(d_1e_k+d_2e_l),$ for some $k\neq l\in\mathbb{N}$ and $\lvert d_1\rvert=\lvert d_2\rvert=1$. In either case, we clearly obtain that $T^{-1}(e_i)$ is not a left-symmetric point, {which establishes} the desired contradiction.
\end{proof}


\begin{thebibliography}{100}
\bibitem{annal}
L. Arambašić, R. Rajić, \textit{``On symmetry of the (strong) Birkhoff–James orthogonality in Hilbert $C^*$-modules"}, \texttt{Ann. Funct. Anal.,
Volume 7, Number 1 (2016), 17-23.}

\bibitem{B} G. Birkhoff,  \textit{``Orthogonality in linear metric spaces"}, \texttt{Duke Math. J., 1 (1935), 169-172.}

\bibitem{CSS} A. Chattopadhyay, D. Sain, T. Senapati, \textit{ ``Characterization of symmetric points in $ l_p^n $-spaces"},\texttt{ Linear Multilinear Algebra, (2019), https://doi.org/10.1080/03081087.2019.1702916.}

\bibitem{Comfort} {W.  W.  Comfort, Ultrafilters : \textit{``Some old  and  some new  results"}, \texttt{Bull. Amer. Math. Soc.,  83 (1977), 417-455.}}

\bibitem{Comfort N} {W. W. Comfort and  S. Negrepontis, \textit{``The theory of ultrafilters"}, \texttt{Grundlehren Math. Wiss., vol. 211, Springer-Verlag, New York, 1974.}}

\bibitem{dkp}
P. Ghosh, D. Sain and K. Paul, \textit{``On symmetry of Birkhoff-James orthogonality of linear operators"}, \texttt{Adv. Oper. Theory, 2 (2017), 428-434.}

\bibitem{1}
 P. Ghosh, K. Paul and D. Sain, \textit{``Symmetric properties of orthogonality of
linear operators on $(\mathbb{R}^
n,\|.\|_1)$"}, \texttt{Novi Sad J. Math., 47 (2017), 41-46.}

\bibitem{james2} R.C. James, \textit{``Inner product in normed linear spaces"}, \textit{Bull. Amer. Math. Soc., 53 (1947), 559-566.}

\bibitem{james} R.C. James, \textit{``Orthogonality and linear functionals in normed linear spaces"}, \texttt{Trans. Amer. Math. Soc., 61 (1947), 265-292.}

\bibitem{3}
N. Komuro, K.-S. Saito and R. Tanaka, \textit{``Left symmetric points for Birkhoff
orthogonality in the preduals of von Neumann algebras"}, \texttt{Bull. Aust. Math.
Soc., 98 (2018), 494-501.}

\bibitem{4}
 N. Komuro, K.-S. Saito and R. Tanaka, \textit{``Symmetric points for (strong)
Birkhoff orthogonality in von Neumann algebras with applications to preserver problems"}, \texttt{J. Math. Anal. Appl., 463 (2018), 1109-1131.}

\bibitem{5}
 N. Komuro, K.-S. Saito and R. Tanaka, \textit{``On symmetry of Birkhoff orthogonality in the positive cones of $C^*$-algebras with applications"}, \texttt{J. Math. Anal.
Appl., 474 (2019), 1488–1497.}

\bibitem{BL}
J. Lamperti, \textit{``On the isometries of certain function-spaces"}, \texttt{Pacific J. Math., 8 (1958), no. 3, 459-466.}

\bibitem{KP} K. Paul, A. Mal and P. W\'{o}jcik, \textit{``Symmetry of Birkhoff-James orthogonality of operators defined between infinite dimensional Banach spaces"}, \texttt{Linear Algebra Appl., {\bf}563 (2019), 142-153.}

\bibitem{Sain2} D. Sain, \textit{``Birkhoff-James orthogonality of linear operators on finite dimensional Banach spaces"}, \texttt{J. Math. Anal. Appl., 447, Issue 2, (2017),  860-866.}

\bibitem{Sain} D. Sain, \textit{``On the norm attainment set of a bounded linear operator"}, \texttt{J. Math. Anal. Appl., 457, Issue 1, (2018), 67-76.}

\bibitem{8}
D. Sain, P. Ghosh and K. Paul, \textit{``On symmetry of Birkhoff-James orthogonality of linear operators on finite-dimensional real Banach spaces"}, \texttt{Oper.
Matrices, 11 (2017), 1087-1095.}

\bibitem{10}
D. Sain, K. Paul, A. Mal, A. Ray, \textit{``A complete characterization of smoothness
in the space of bounded linear operators"}, \texttt{Linear Multilinear Algebra, (2019),
doi.org/10.1080/03081087.2019.1586824.}

\bibitem{SRBB} D. Sain, S. Roy, S. Bagchi and V. Balestro, \textit{``A study of symmetric points in Banach spaces"}, \texttt{Linear Multilinear Algebra, (2020), https://doi.org/10.1080/03081087.2020.1749541.}

\bibitem{12}
A. Turn\^sek, \textit{``A remark on orthogonality and symmetry of operators in B(H)"},
\texttt{Linear Algebra Appl., 535 (2017), 141-150.}

\bibitem{turnsek}
 A. Turnsek, \textit{``On operators preserving James’ orthogonality"}, \texttt{Linear Algebra and its Applications, 407 (2005), 189-195.}
\end{thebibliography}


\end{document}